\documentclass[conference]{IEEEtran}
\IEEEoverridecommandlockouts
\usepackage{cite}
\usepackage{amsmath,amssymb,amsfonts}
\usepackage{algorithmic}
\usepackage{graphicx}
\usepackage{textcomp}
\usepackage{xcolor}
\def\BibTeX{{\rm B\kern-.05em{\sc i\kern-.025em b}\kern-.08em
  T\kern-.1667em\lower.7ex\hbox{E}\kern-.125emX}}

\usepackage{cmap}
\usepackage[T1,T2A]{fontenc}
\usepackage[utf8]{inputenc}

\begin{document}

\title{Formal asymptotics of parametric subresonance}

\author{\IEEEauthorblockN{1\textsuperscript{st}Polina Astafyeva}
\IEEEauthorblockA{\textit{Inst. of Math. with Com. Centre of UFRC RAS,} \\
\textit{Ufa State Petroleum Technological University}\\
Ufa, Russia \\
astafyeva.polina2@gmail.com}
\and
\IEEEauthorblockN{2\textsuperscript{nd} Oleg Kiselev}
\IEEEauthorblockA{\textit{Innopolis University,} \\
\textit{Inst. of Math. with Com. Centre of UFRC RAS}\\
Innopolis, Russia \\
o.kiselev@innopolis.ru}

}

\maketitle

%\IEEEpubidadjcol

%\begin{abstract}
%The behavior of a linear oscillator under the action of an external almost periodic force is investigated. The constructed solutions grow more slowly than the resonant ones. The dependence of the amplitude of growing solutions on the parameters of an almost periodic perturbation is calculated.
%\end{abstract}

%\begin{IEEEkeywords}
%Classical Analysis and ODEs, subresonant
%\end{IEEEkeywords}

%\section{Introduction}

\label{secIntroduection}
Lots of kinds of linear equations of second order are the most important models of mechanics because of Newtonian equations of a motion. Among such equations we are interested in an equation with time-dependent coefficient. In this work we consider two different equations of such kind. The first one is the equation with almost periodic coefficient:  
\begin{equation}
u''+( \omega^2 +\epsilon q(t) ) u=0 
\label{iskh}
\end{equation}
here $q(t)$ is almost periodic function and $\epsilon$ is a small positive parameter.

This equation has two important properties which define behaviours of solutions. The first one is a coefficient $\omega^2$. It defines an oscillation of the solution for the simplest case $q(t)\equiv0$. 

Another coefficient $q(t)$ can or cannot imply the resonant behaviour for the solution.  The key work for understanding the behaviour of the solution of such kind of equation is the Floquet theory \cite{Floquet1883}. The parametric map for the parameters $(\omega,\epsilon)$ of the equation with periodic coefficient is called Arnold's tonques \cite{Arnold1963}. Such map defines the zones of parametric resonances \cite{BogolyubovMitropolskii1961Eng}. 

One of the famous examples of such equation is Mathieu equation for special kind of the coefficient:
$$
q(t)\equiv a\cos(2t).
$$
Although the Mathieu equation appeared in theory of Laplace equation with elliptic boundary \cite{WhittakerWatson1920} now this equation is widely used model in quantum mechanics and theory of parametric resonance \cite{BogolyubovMitropolskii1961Eng},\cite{Sternkopf}, \cite{def},\cite{orb}.

One of more general equations of such kind appeared in the theory of Moon motion \cite{Hill}, where $q(t)$ is a periodic function defined by their Fourier series.

Here we study more general and more natural case when the $q(t)$ is almost-periodic. Such sight on the coefficient looks as more physically important  since in real phenomenon pure periodic coefficient is enveloped by additional disturbances and perturbations. Moreover, the period of the coefficient often does not coincide with period of  frequency of solution for the equation, such as $q(t)\equiv0$.

Let us consider an second order equation with almost periodic coefficient:
\begin{equation}
u''+( \omega^2 +\epsilon q(t) ) u=0 
\label{iskh}
\end{equation}
Here almost periodic function \cite{Levitan1953}, \cite{LevitanZhikov1977}:
\begin{equation}
q(t)=\sum_{n=1}^{\infty}\frac{1}{n^k}\cos((2-\frac{1}{n^p})t).
\label{formulaForF}
\end{equation}
$$k>1, p>0$$

The parameter $\omega$ differs little from 1: $\omega=1+\delta$, here $\delta$-- is the parameter of the equation.
The task is to determine the areas of stability of solutions of the equation depending
on the parameters $\delta$ и $\epsilon$.

In the following section we formulate aims of this work. Next section contains an asymptotic approach to the result. The last section shows the closeness of asymptotic and numerical approaches and contains the conclusions.
  
\section{Statement of the problem}
For the equation with almost-periodic coefficient one can see three different behaviour with respect to parameter $t\to\infty$.

\begin{figure}
 \includegraphics[scale=0.12]{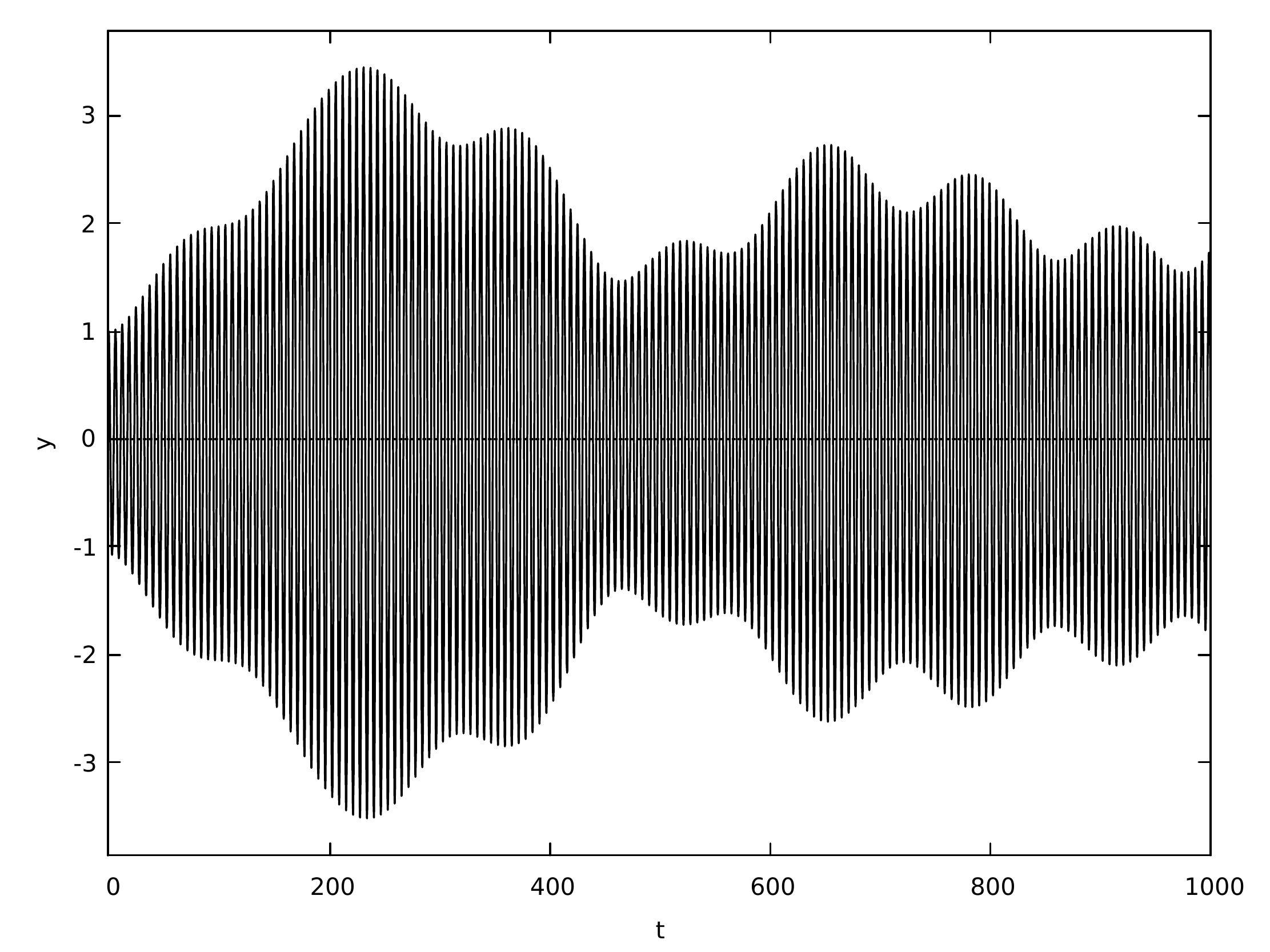}   \includegraphics[scale=0.12]{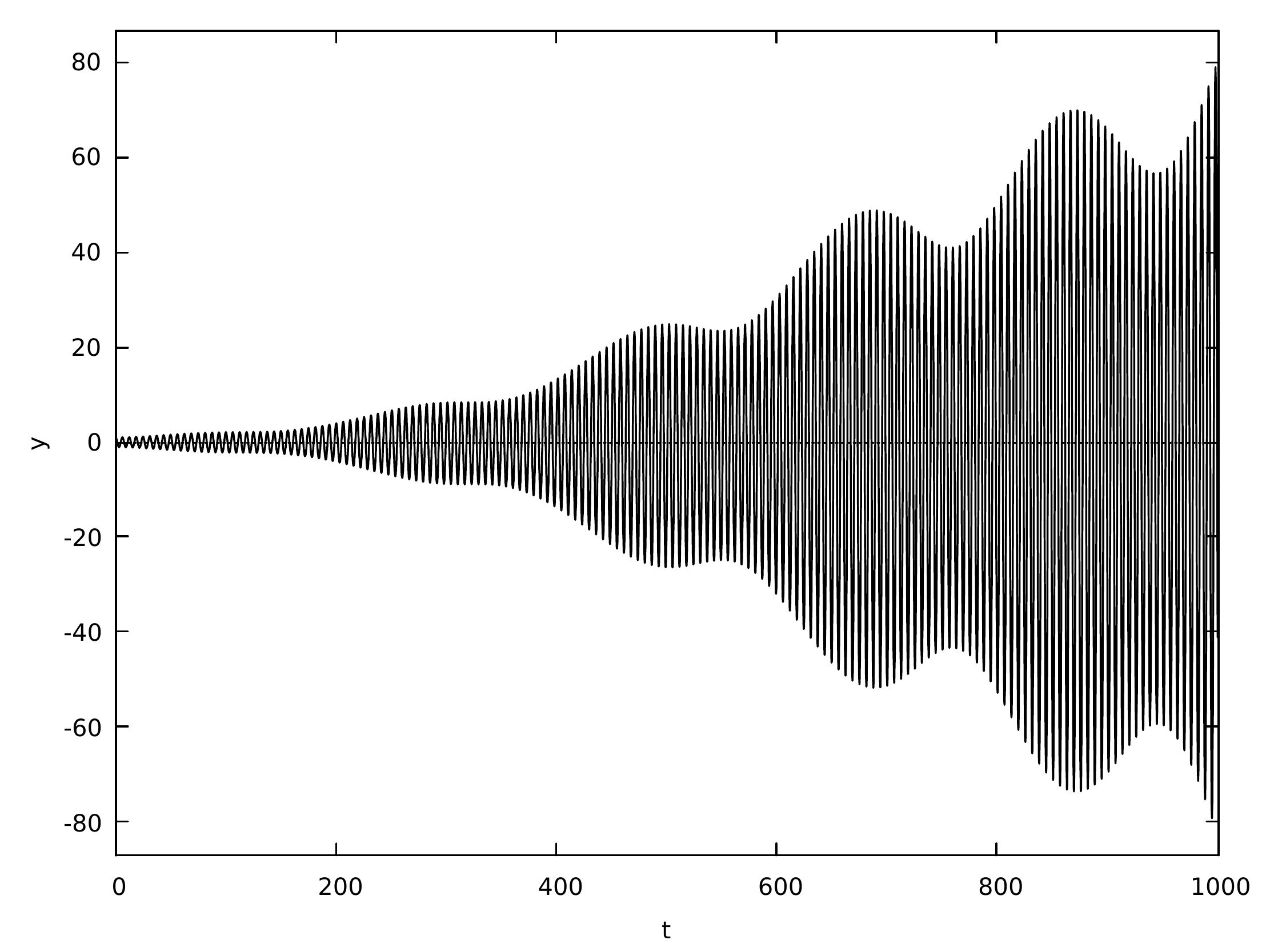}
 \includegraphics[scale=0.12]{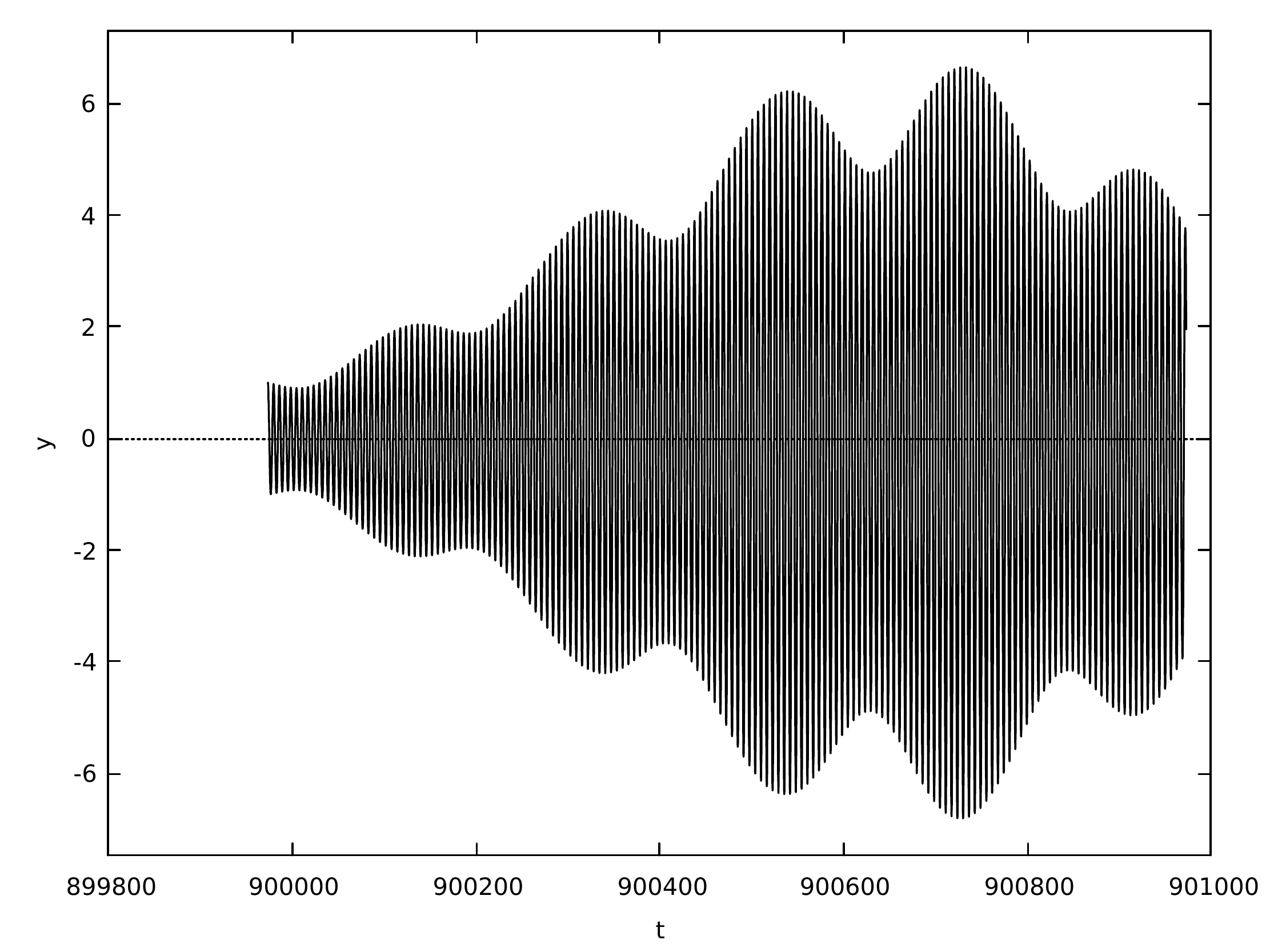}
\caption{Here we give three different cases for the solution of equation (\ref{iskh}) where $k=2$,and $p=5$. The left picture shows bounded oscillations for the equation $\omega=1+\delta$, where $\delta=0.01$, $\epsilon=0.1$. The middle picture shows the growing solution  as $\delta=0.00347$, $\epsilon=0.1$, and the right picture indicate the final stage of this growing solution as the same $\delta=0.00347$, $\epsilon=0.1$  at very large time interval around $t=900500$.}
\end{figure}

Our first problem for studying is to understand a behaviour of solution depending on the coefficients. $\omega$ and$\epsilon$ for given form $q(t)$ as the almost-periodic function. 

Second and more intrigue aim for us is to find an dependency of amplitude of the oscillated solution which one can consider as a similarity to parametric resonance growth which was derived for example in \cite{BogolyubovMitropolskii1961Eng}.

\section{Asymptotic aproach}
Let us construct an asymptotic solution in the form:
\begin{equation}
u \sim u_{0}+\epsilon u_{1},
\label{as}
\end{equation}

Let us substitute (\ref{as}) in (\ref{iskh}) and combine the terms at the same degrees $\epsilon$. Obtain the equation for the main term:
$$\frac{d^2}{dt^2}u_{0}+\omega^2 u_{0}=0 .$$
Let us look for its solution in the form of $u_{0}=a(\tau) \cos(\omega t)+ b(\tau) \sin(\omega t)$, using two scale method \cite{BogolyubovMitropolskii1961Eng}, where $\tau=\epsilon^{\gamma}t$ is slow time. 
The parameter $\gamma >0$,  will be determined below from the condition of uniformity of the asymptotic solution for large values$t$.

The equation for the first correction:
\begin{eqnarray}
\frac{{{d}^{2}}}{d {{t}^{2}}} u_{1}+u_{1}\, {{\omega }^{2}}+q(t)  \left( b \sin{\left( t \omega \right) }+a \cos{\left( t \omega \right) }\right) -
\nonumber\\ -2 \mathit{a_{1}} \omega  \sin{\left( t \omega \right) }+2 \mathit{b_{1}} \omega  \cos{\left( t \omega \right) }=0
\end{eqnarray}
here $a_{1}=\epsilon^{\gamma-1}a'$ and $b_{1}=\epsilon^{\gamma-1}b'$, a stroke means a slow-time derivative $\tau$.
 %(x*y)''=(x'*y+x*y')'=x''*y+2x'*y'*x*y''. 
%$$x=a(\epsilon t), y=\cos(\omega t)$$

The task is to find the dependence on  $\tau$ for a and b.
A limited solution for the first correction of the function $u_1$ s constructed below.. 

Denote
$$
f(t)=-q(t) \left( b \sin{\left( t \omega \right) }+a \cos{\left( t \omega \right) }\right) 
$$
   In this case, the general solution for the first correction equation can be written as:

\begin{eqnarray}u_{1}=A\cos(\omega t)+B\sin(\omega t)-b_{1} t \sin(t \omega)-a_{1} t \cos(t \omega)+
\nonumber\\
+ \frac{\cos(\omega t)}{\omega} \int_{0}^{t}{f(\tilde{t})\sin(\omega \tilde{t})} d\tilde{t} -\frac{\sin(\omega t)}{\omega}\int_{0}^{t}{f(\tilde{t})\cos(\omega \tilde{t})} d\tilde{t}.
\label{ob}
\end{eqnarray}

let us introduce a replacement in the integrand function $\omega \tilde{t}\equiv \tilde{t}+\delta \tilde{t} =\tilde{t}+\kappa \tau$

  Note that the formula for  $f (\tilde{t})$ contains a sum over n, we change the summation with integration over $\tilde{t}$.
Then the solution of the equation for  $u_1$ is represented as the sum of n for the first integral in the formula (\ref{ob}): 

\begin{eqnarray*}\frac{a \sin{\left( \frac{2 {{n}^{p}} \kappa  \tau +\left( 4 {{n}^{p}}-1\right)  \tilde{t}}{{{n}^{p}}}\right) }}{4 {{n}^{k}}}-\frac{b \cos{\left( \frac{2 {{n}^{p}} \kappa  \tau +\left( 4 {{n}^{p}}-1\right)  \tilde{t}}{{{n}^{p}}}\right) }}{4 {{n}^{k}}}+\nonumber\\+
\frac{a \sin{\left( \frac{2 {{n}^{p}} \kappa  \tau +\tilde{t}}{{{n}^{p}}}\right) }}{4 {{n}^{k}}}-\frac{b \cos{\left( \frac{2 {{n}^{p}} \kappa  \tau +\tilde{t}}{{{n}^{p}}}\right) }}{4 {{n}^{k}}}+\frac{b \cos{\left( \frac{\left( 2 {{n}^{p}}-1\right)  \tilde{t}}{{{n}^{p}}}\right) }}{2 {{n}^{k}}}
\end{eqnarray*}

and the second integral in the formula (\ref{ob}):
\begin{eqnarray*}
\frac{b \sin{\left( \frac{2 {{n}^{p}} \kappa  \tau +\left( 4 {{n}^{p}}-1\right)  \tilde{t}}{{{n}^{p}}}\right) }}{4 {{n}^{k}}}+\frac{a \cos{\left( \frac{2 {{n}^{p}} \kappa  \tau +\left( 4 {{n}^{p}}-1\right)  \tilde{t}}{{{n}^{p}}}\right) }}{4 {{n}^{k}}}+
\nonumber\\+
\frac{b \sin{\left( \frac{2 {{n}^{p}} \kappa  \tau +\tilde{t}}{{{n}^{p}}}\right) }}{4 {{n}^{k}}}+\frac{a \cos{\left( \frac{2 {{n}^{p}} \kappa  \tau +\tilde{t}}{{{n}^{p}}}\right) }}{4 {{n}^{k}}}+\frac{a \cos{\left( \frac{\left( 2 {{n}^{p}}-1\right)  \tilde{t}}{{{n}^{p}}}\right) }}{2 {{n}^{k}}}
\end{eqnarray*}

After  sequential integration  of all parts by  $\tilde{t}$ let us allocate the maximum order  with respect to $n$ for obtained formula :

\begin{eqnarray*}
\frac{\cos(\omega t)}{{4 \omega}} \sum_{n=1}^{\infty} n^{p-k}\left(b \left( \sin{\left( 2 \kappa  \tau \right) }- \sin{\left( 2 \kappa  \tau +\frac{t}{{{n}^{p}}}\right) } \right)+\right.
\\
\left. a \left( \cos{\left( 2 \kappa  \tau \right) }- \cos{\left( 2 \kappa  \tau +\frac{t}{{{n}^{p}}}\right) }  \right) \right) 
\end{eqnarray*}

and 
\begin{eqnarray*}
\frac{\sin(\omega t)}{{4 \omega}} \sum_{n=1}^{\infty} n^{p-k}\left(a \left( \sin{\left( 2 \kappa  \tau \right) }- \sin{\left( 2 \kappa  \tau +\frac{t}{{{n}^{p}}}\right) } \right)-\right.
\\
\left.b \left( \cos{\left( 2 \kappa  \tau \right) }- \cos{\left( 2 \kappa  \tau +\frac{t}{{{n}^{p}}}\right) }  \right) \right) 
\end{eqnarray*}

%К растущим слагаемым приводят интегралы от :

%\[\frac{a \sin{\left( \frac{2 {{n}^{p}} \kappa  \tau +t}{{{n}^{p}}}\right) }}{4 {{n}^{k}}}-\frac{b \cos{\left( \frac{2 {{n}^{p}} \kappa  \tau +t}{{{n}^{p}}}\right) }}{4 {{n}^{k}}}+\frac{b \sin{\left( \frac{2 {{n}^{p}} \kappa  \tau +t}{{{n}^{p}}}\right) }}{4 {{n}^{k}}}+\frac{a \cos{\left( \frac{2 {{n}^{p}} \kappa  \tau +t}{{{n}^{p}}}\right) }}{4 {{n}^{k}}}\]

%преобразуем:
%\[-\frac{\left( a\, {{n}^{p}} \sin{\left( \frac{2 {{n}^{p}} \kappa  \tau +t}{{{n}^{p}}}\right) }-b\, {{n}^{p}} \cos{\left( \frac{2 {{n}^{p}} \kappa  \tau +t}{{{n}^{p}}}\right) }-a\, {{n}^{p}} \sin{\left( 2 \kappa  \tau \right) }+b\, {{n}^{p}} \cos{\left( 2 \kappa  \tau \right) }\right)  \sin{\left( t \omega \right) }}{4 {{n}^{k}} \omega }-\frac{\left( b\, {{n}^{p}} \sin{\left( \frac{2 {{n}^{p}} \kappa  \tau +t}{{{n}^{p}}}\right) }+a\, {{n}^{p}} \cos{\left( \frac{2 {{n}^{p}} \kappa  \tau +t}{{{n}^{p}}}\right) }-b\, {{n}^{p}} \sin{\left( 2 \kappa  \tau \right) }-a\, {{n}^{p}} \cos{\left( 2 \kappa  \tau \right) }\right)  \cos{\left( t \omega \right) }}{4 {{n}^{k}} \omega }\]

In the solution  $u_1$ In the solution $u_1$ the growing terms in t are combined at $\cos(\omega t)$:
\begin{eqnarray*}-\frac{a}{4 \omega}\left(\cos{\left( 2 \kappa  \tau \right) } \sum_{n=1}^{\infty} {{n}^{p-k}}\,  \left( \cos{\left( \frac{t}{{{n}^{p}}}\right) }-1\right) -\right.
\\
\left.
\sin{\left( 2 \kappa  \tau \right) }\sum_{n=1}^{\infty}\sin{\left( \frac{t}{{{n}^{p}}}\right) }\right)  +
\nonumber\\-\frac{b}{4 \omega}\left(\sin{\left( 2 \kappa  \tau \right) }\sum_{n=1}^{\infty} {{n}^{p-k}}\,  \left( \cos{\left( \frac{t}{{{n}^{p}}}\right) }-1\right)  +\right.
\\
\left.\cos{\left( 2 \kappa  \tau \right)}\sum_{n=1}^{\infty}\sin{\left( \frac{t}{{{n}^{p}}}\right) }  \right) -\mathit{a_{1}} t
\end{eqnarray*}

and $\sin(\omega t)$:

\begin{eqnarray*}
\frac{b}{4 \omega}\left(\sin{\left( 2 \kappa  \tau \right) }\sum_{n=1}^{\infty}{{n}^{p-k}}\,  \sin{\left( \frac{t}{{{n}^{p}}}\right) } +\right.
\\
\left.\cos{\left( 2 \kappa  \tau \right)\sum_{n=1}^{\infty}\left( 1-\cos{\left( \frac{t}{{{n}^{p}}}\right) }\right)   }\right) +
\nonumber\\-\frac{a}{4 \omega}\sin{\left( 2 \kappa  \tau \right) }\left( \sum_{n=1}^{\infty}{{n}^{p-k}}\, \left( \cos{\left( \frac{t}{{{n}^{p}}}\right) }-1\right)  -\right.
\\
\left.\cos{\left( 2 \kappa  \tau \right) }\sum_{n=1}^{\infty}\sin{\left( \frac{t}{{{n}^{p}}}\right) } \right) -\mathit{b_{1}} t
\end{eqnarray*}

Replace: 

$$\sum_{n=1}^{\infty}{{n}^{p-k}}\, \left( \cos{\left( \frac{t}{{{n}^{p}}}\right) }-1\right)=\sum_{n=1}^{\infty}n^{p-k} 2\sin \left( \frac{t}{2n^p}\right)$$

In \cite{AstafyevaKiselev}, the asymptotics were calculated:
\begin{eqnarray*}
\sum_{n=1}^{\infty}n^{p-k} \sin(\frac{t}{n^p})\sim t^{1-\alpha} C_s, 
\\
\qquad C_s\sim\frac{1}{p} 
\int_{0}^{\frac{\pi }{2}}\tau^{\alpha-2} \sin(\tau)d\tau,
\\
\sum_{n=1}^{\infty} n^{p-k} 2\sin^2(\frac{t}{2n^p})\sim t^{1-\alpha} C_c, 
\\
\qquad C_c\sim-\left(\frac{1}{2\pi^{1-\alpha}p(1-\alpha)} + 
\frac{1}{p}\int_0^{\pi}\tau^{\alpha-2} \sin^2( \tau/2)d\tau\right). 
\end{eqnarray*}

Denote:
$$\alpha=\frac{k-1}{p}.$$
Let us combine the growing terms in t at$\cos(\omega t)$ и $\sin(\omega t)$, in order for the coefficients with increasing summands to be zeros:

\begin{eqnarray}\epsilon^{\gamma-1}ta'=\frac{a t^{1-\alpha}}{4 \omega }\, \left( \mathit{C_c} \cos{\left( 2 \kappa  \tau \right) }-\mathit{C_{s}} \sin{\left( 2 \kappa  \tau \right) }\right) +
\nonumber\\+\frac{b t^{1-\alpha}}{4\omega }\, \left( \mathit{C_c} \sin{\left( 2 \kappa  \tau \right) }+\mathit{C_s} \cos{\left( 2 \kappa  \tau \right) }\right) 
\end{eqnarray}

\begin{eqnarray}\epsilon^{\gamma-1}tb'=\frac{bt^{1-\alpha}}{4 \omega }\, \left( \mathit{C_{s}} \sin{\left( 2 \kappa  \tau \right) }-\mathit{C_c} \cos{\left( 2 \kappa  \tau \right) }\right) +
\nonumber\\+\frac{a t^{1-\alpha}}{4\omega }\, \left( \mathit{C_c} \sin{\left( 2 \kappa  \tau \right) }+\mathit{C_s} \cos{\left( 2 \kappa  \tau \right) }\right) 
\end{eqnarray}

In these formulas, we will reduce both parts by t and make a replacement $t=\frac{\tau}{\epsilon^{\gamma}}$ .
As a result , we get:

\begin{eqnarray*}\epsilon^{\gamma-1}\frac{d}{d \tau } a=\frac{a \epsilon^{\gamma \alpha}}{4 {{\tau }^{\alpha }} \omega }\, \left( \mathit{C_c} \cos{\left( 2 \kappa  \tau \right) }-\mathit{C_{s}} \sin{\left( 2 \kappa  \tau \right) }\right) +\nonumber\\+\frac{b\epsilon^{\gamma \alpha}}{4 {{\tau }^{\alpha }} \omega }\, \left( \mathit{C_c} \sin{\left( 2 \kappa  \tau \right) }+\mathit{C_s} \cos{\left( 2 \kappa  \tau \right) }\right),\\
\epsilon^{\gamma-1}\frac{d}{d \tau } b=\frac{b \epsilon^{\gamma \alpha}}{4 {{\tau }^{\alpha }} \omega }\, \left( \mathit{C_s} \sin{\left( 2 \kappa  \tau \right) }-\mathit{C_c} \cos{\left( 2 \kappa  \tau \right) }\right) +\nonumber\\+\frac{a\epsilon^{\gamma \alpha}}{4 {{\tau }^{\alpha }} \omega }\, \left( \mathit{C_c} \sin{\left( 2 \kappa  \tau \right) }+\mathit{C_s} \cos{\left( 2 \kappa  \tau \right) }\right) .
\end{eqnarray*}
Hence $$\gamma-1=\gamma \alpha, \quad \gamma=\frac{1}{1-\alpha}.$$ 
As a result , the system of equations will take the form:
\begin{eqnarray*}\frac{d}{d \tau } a=\frac{a}{4 {{\tau }^{\alpha }} \omega }\, \left( \mathit{C_c} \cos{\left( 2 \kappa  \tau \right) -\mathit{C_{s}} \sin{\left( 2 \kappa  \tau \right) }}\right) +
\nonumber\\+
\frac{b}{4 {{\tau }^{\alpha }} \omega }\, \left( \mathit{C_c} \sin{\left( 2 \kappa  \tau \right) }+\mathit{C_s} \cos{\left( 2 \kappa  \tau \right) }\right), \\
\frac{d}{d \tau } b=\frac{b}{4 {{\tau }^{\alpha }} \omega }\, \left( \mathit{C_s} \sin{\left( 2 \kappa  \tau \right) }-\mathit{C_c} \cos{\left( 2 \kappa  \tau \right) }\right) +
\nonumber\\+\frac{a}{4 {{\tau }^{\alpha }} \omega }\, \left( \mathit{C_c} \sin{\left( 2 \kappa  \tau \right) }+\mathit{C_s} \cos{\left( 2 \kappa  \tau \right) }\right). 
\end{eqnarray*}

%Для того чтобы порядки в правой и левой частей уравнений совпадали, примем:
%$$\gamma=\frac{1}{1-\alpha}.$$
Denote:
$$
A_{\alpha}=\sqrt{C_s^2+C_c^2},\quad 
\phi_{\alpha}=\arctan\left(\frac{C_c}{C_s}\right),
$$
Then the system will take the form:

\[\frac{d}{d \tau } a=\frac{A\, \left( b \sin{\left(  \phi +2 \kappa  \tau \right) }+a \cos{\left( \phi +2 \kappa  \tau \right) }\right) }{4 \omega
{{\tau }^{\alpha }}}\]

\[\frac{d}{d \tau } b=\frac{A\, \left( a \sin{\left(  \phi +2 \kappa  \tau \right) }-b \cos{\left(  \phi +2 \kappa  \tau \right) }\right) }{4 \omega{{\tau }^{\alpha }}}\]

Denote $\frac{A}{4\omega}=B$, then:

\[\frac{d}{d \tau } a=\frac{B\, \left( b \sin{\left(  \phi +2 \kappa  \tau \right) }+a \cos{\left( \phi +2 \kappa  \tau \right) }\right) }{
{{\tau }^{\alpha }}}\]

\[\frac{d}{d \tau } b=\frac{B\, \left( a \sin{\left(  \phi +2 \kappa  \tau \right) }-b \cos{\left(  \phi +2 \kappa  \tau \right) }\right) }{{{\tau }^{\alpha }}}\]

The system can also be written in a complex form:
\[\frac{d}{d \tau } z=\frac{B\, \overline{z}\, {{ e}^{i \phi +2 i \kappa  \tau }}}{{{\tau }^{\alpha }}}\]

$\delta=\epsilon^{\gamma}\kappa$

Denote $z=Ze^{i\kappa\tau+\frac{i\phi}{2}}$

\[i Z \kappa \, {{e}^{\frac{ i \phi }{2}+i \kappa  \tau }}+\left( \frac{d}{d \tau } Z\right) \, {{ e}^{\frac{ i \phi }{2}+ i \kappa  \tau }}=\frac{B\, \overline{Z}\, {{ e}^{\frac{ i \phi }{2}+i \kappa  \tau }}}{{{\tau }^{\alpha }}}\]

simplify:
\[i Z \kappa +\frac{d}{d \tau } Z=\frac{B\, \overline{Z}}{{{\tau }^{\alpha }}}\]

\[\frac{d}{d \tau } Z=\frac{B\, \overline{Z}}{{{\tau }^{\alpha }}}- i Z \kappa \]

Let us write Z as the sum of the real and imaginary parts:  $Z=w+iv.$ 

As a result , a system of equations is obtained:
\begin{eqnarray}
\frac{d}{d \tau } w=\frac{B w}{{{\tau }^{\alpha }}}+v\, \kappa,
\nonumber\\
\frac{d}{d \tau } v=-w \kappa-\frac{B v}{{{\tau }^{\alpha }}}.
\label{syst}
 \end{eqnarray}

\section{Comparison of the asymptotic approaches and numerical results}
Let us consider the properties of derived system (\ref{syst})
for different values of the parameters. 

In the case $\kappa=0$  the system splits into two equations. Their solutions have the following form:

\begin{eqnarray*}
w=C_1\exp\left (\frac{B\tau^{1-\alpha}}{1-\alpha}\right), \quad v=C_2\exp\left (\frac{-B\tau^{1-\alpha}}{1-\alpha}\right),
\end{eqnarray*}
Here one can see the solution is exponentially growing because of $w(t)$.

For$\kappa \not=0$ we divide both parts of the equations on the parameter  $\kappa$ and rewrite $\kappa\tau=\theta$. In this case we get:
\begin{eqnarray*}
\frac{d}{d \theta} w=\frac{B w}{\kappa^{1-\alpha}\,\theta ^{\alpha }}+v\,  \\
\frac{d}{d \theta } v=-\frac{B v}{\kappa^{1-\alpha}\,\theta ^{\alpha }}-w 
\end{eqnarray*}
Let us rewrite $B\kappa^{\alpha-1}=\lambda$, then:
\begin{eqnarray*}
\frac{d}{d \theta} w=\frac{\lambda}{\theta ^{\alpha }}w+v,
\frac{d}{d \theta } v=-\frac{\lambda}{\theta ^{\alpha }}v-w.
\end{eqnarray*}
Here $\lambda>0$ is a parameter of the equation and  $\kappa \to 0, \,\lambda \to \infty$.

This system leads to the second order differential equation:
\begin{eqnarray*}
\frac{d^2}{d \theta^2} w+ \left(1- \frac{\lambda^2}{\theta ^{2\alpha }}+\frac{\alpha \lambda}{\theta ^{\alpha+1 }}\right)w=0
\end{eqnarray*}
Here $\lambda$ is a large parameter. So asymptotic solution can be obtained by the WKB method \cite{wkb}:
$$
w\sim C_1\frac{\exp(\int^t\sqrt{1- \frac{\lambda^2}{\theta ^{2\alpha }}+\frac{\alpha \lambda}{\theta ^{\alpha+1 }}} dt)}{\sqrt[4]{1- \frac{\lambda^2}{\theta ^{2\alpha }}+\frac{\alpha \lambda}{\theta ^{\alpha+1 }}}}.
$$ 
So the turning point to change the growing character of the solution is located in a neighbourhood of the point 
$$
\theta*\sim\lambda^{1/\alpha}.
$$

Let us compare the position of the tuning point obtained here by asymptotical   way and the results which ware shown in the figure. 

For the numerical example we have following data 
$$
k=2,\quad p=5,\quad \alpha\equiv(k-1)/p=1/5,\quad \gamma\equiv 1/(1+\alpha)=6/5,
$$
$$
\quad A\sim1.09264275,
$$
here the value of $A$ was obtained with respect to the work \cite{AstafyevaKiselev}.
$$
B\equiv\frac{A}{4\omega}\sim 1/4,\quad \lambda=5, 
\kappa\equiv(B/\lambda)^(1/(1-\alpha))\sim \left(\frac{1}{4}\right)^{5/4},
$$
$$
\epsilon=0.1,
\quad \delta=\epsilon^{\gamma\kappa},
\quad \theta_*\sim\lambda^{1/\alpha},
\quad 
T=\frac{\theta_*}{\kappa\epsilon^\gamma}\sim 900473.
$$ 
 
{\bf Conclusions.} We obtain the system of equations for parametric sub-resonant growth of the amplitude of oscillations. The growth depends on the slow variable $\tau=\epsilon^{\gamma} t$. We find also the time of turning point form the growing of the amplitude to the bounded oscillations in the slow variable $\tau$. The comparison between the asymptotic approximation for the turning time and numerical one is shown.


\begin{thebibliography}{00}
\bibitem{Floquet1883}
Floquet, ``Sur les equations differentielles lineaires a cofficients
	periodiques``, Ann. de l'Ecole norm. sup., 1883, XII, pp. 47-88.
 
\bibitem{Levitan1953} Levitan, `` Almost periodic functions,'' Gos. izd. tech.-teor. lit., 1953.
\bibitem{LevitanZhikov1977} B.M.Levitan, V.V.Zhikov.  ``Almost periodic function and differential equations,''  M. MSU, 1978.
\bibitem{Arnold1963}V. I. Arnold , ``Small Denominators and problems of stability of motion in classical and celestial mechanics,`` Uspekhi Mat. Nauk, 18(6):91--192, 1963.

\bibitem{BogolyubovMitropolskii1961Eng}N.N. Bogolyubov and Yu.A. Mitropolskii, ``Asymptotic methods in the theory of non-linear oscillations``, Gordon and Breach science publishers, New York,1961, p. 537.
\bibitem{Sternkopf}Sternkopf, C.; Diethold, C.; Gerhardt, U.; Wurmus, J., Manske, E. ``Heterodyne interferometer laser source with a pair of two phase locked loop coupled HetextendashNe lasers by 632.8 nm `` Measurement Science and Technology, IOP Publishing, 2012, 23, 074006

\bibitem{WhittakerWatson1920}E. T. Whittaker and G. N. Watson ``A course of modern analysis``, Cambridge University Press, 1920.

\bibitem{Hill}G.W. Hill, ``On the part of the motion of the lunar perigee which is a function of the mean motions of the sun
and moon``, Cambridge, Mass., Press of J. Wilson and son, 1877, Cambridge, Mass., Press of J. Wilson and son. 
\bibitem{AstafyevaKiselev}P.Y.Astafyeva, O.M.Kiselev, ``Subresonant solutions of the linear oscillator equation, 2021 International Conference Nonlinearity, Information and Robotics, Innopolis, Russia, pp 90-94.
\bibitem{def}Artem Eremin, Mikhail Golub, Evgeny Glushkov, Natalia Glushkova.  ``Identification of delamination based on the Lamb wave scattering resonance frequencies`` NDT and E International Volume 103, April 2019, Pages 145-153
\bibitem{orb} Peale, S. J. ``Orbital Resonances in the Solar System`` Annual Review of Astronomy and Astrophysics Vol. 14:215-246 

\bibitem{wkb}Olver, Frank William John , Asymptotics and Special Functions,1974, Academic Press.
\end{thebibliography}
\end{document}